\documentclass[11pt]{article}
\usepackage{a4wide}
\usepackage{amsfonts}
\usepackage{amsmath}
\usepackage{amssymb}
\usepackage{latexsym}
\usepackage{amsthm}
\def\newf{\kern-.15em{f}}
\def\newg{\kern-.15em{g}}
\def\newj{\kern-.15em{j}}
\def\lceiling{\lceil}
\def\rceiling{\rceil}
\def\Q{\mathbb Q}
\def\R{\mathbb R}
\def\W{\mathcal W}
\def\F{\mathcal F}
\def\Ot{\mathcal O}
\def\id{\mathrm{Id}}
\def\E{\mathcal E}
\def\Ec{\E^{\mathrm{cusp}}}
\def\Ee{\E^{\mathrm{Eis}}}

\def\Id{\id}

\def\OC2{\O_{\C_2}}

\def\<{\langle}
\def\>{\rangle}

\def\Z{\mathbb Z}
\def\T{\mathbf T}
\def\N{\mathbb N}
\def\A{\mathbb A}
\def\C{\mathbb C}
\def\X{\mathcal X}
\DeclareMathOperator{\Sp}{Sp}
\DeclareMathOperator{\Hom}{Hom}
\DeclareMathOperator{\GL}{GL}
\DeclareMathOperator{\pr}{pr}
\title{The $2$-adic Eigencurve is  Proper.}
\author{Kevin Buzzard \and Frank Calegari\footnote{Supported
in part by the American Institute of Mathematics}}
\begin{document}
\maketitle

%%%%%%% Introduction

\newtheorem{theorem}{Theorem}[section]
\newtheorem{df}[theorem]{Definition}
\newtheorem{lemma}[theorem]{Lemma}
\newtheorem{prop}[theorem]{Proposition}
\newtheorem{remark}[theorem]{Remark}
\newtheorem{cor}[theorem]{Corollary}

\section{Introduction}

In~\cite{Eigencurve}, Coleman and Mazur construct a rigid analytic
space $\E$ that parameterizes overconvergent and therefore classical
modular eigenforms of finite slope.
The geometry of $\E$ is at present poorly understood, and seems quite
complicated, especially over the centre of weight space. Recently,
some progress has been made in understanding the geometry of $\E$ 
in certain examples (see for example~\cite{buzzcal},\cite{buzzkill}).
Many questions remain. In this paper, we address the following question
raised on p5 of~\cite{Eigencurve}:

\begin{quote}
Do there exist $p$-adic analytic families of overconvergent eigenforms of
finite slope parameterized by a punctured disc, and converging, at the 
puncture, to an overconvergent eigenform of infinite slope?
\end{quote}

We answer this question in the negative for the 2-adic eigencurve of tame
level~1. Another
way of phrasing our result is that the map from the eigencurve to weight
space  satisfies the \emph{valuative criterion of properness}, and it
is in this sense that the phrase  
``proper'' is used in the title, since the
projection to weight space has infinite degree and so
 is not technically  proper in the
sense of rigid analytic geometry.
One might perhaps say that this map is ``functorially
proper''.
Our approach is based on the following simple idea. One knows 
(for instance, from \cite{wild}) that finite
slope eigenforms of integer weight may be analytically continued
far into the supersingular regions of the moduli space.
On the other hand, it turns out that eigenforms
in the kernel of $U$ do not extend as far.
Now one can check that a limit of highly overconvergent eigenforms is also
highly overconvergent, and this shows that the given a punctured disc as
above, the limiting eigenform cannot lie in the kernel of $U$.

The problem with this approach is that perhaps the most natural definition
of ``highly convergent'' is not so easy to work with at non-integral
weight. The problem stems from
the fact that such forms of non-integral weight are not defined as sections of
a line bundle. In fact Coleman's definition of an overconvergent
form of weight $\kappa$ is a formal $q$-expansion $F$ for
which $F/E_{\kappa}$ is overconvergent of weight~$0$, where $E_\kappa$
is the weight $\kappa$ $p$-deprived Eisenstein series. One might then
hope that the overconvergence of $F/E_{\kappa}$ would be a good measure
of the overconvergence of $F$. One difficulty is that
if~$F$ is an eigenform for the Hecke operators,
the form $F/E_{\kappa}$ is unlikely to
be an eigenform. This does not cause too much trouble when proving
that finite slope eigenforms overconverge a long way, as one can
twist the $U$-operator as explained in~\cite{Coleman} and apply
the usual techniques. We outline the argument in sections 2 and 3
of this paper. On the other hand we do not know how to prove
general results about (the lack of) overconvergence of forms
in the kernel of $U$ in this generality. Things would be easier
if we used $V(E_\kappa)$ to twist from weight~$\kappa$ to
weight~0, but unfortunately the results we achieve using
this twist are not strong enough for us to get the strict
inequalities that we need.

The approach that we take in our ``test case'' of $N=1$ and $p=2$
is to control the kernel of $U$ in weight $\kappa$
by explicitly writing down the matrix
of $U$ (and of $2VU-\Id$) with respect to a carefully-chosen basis.
To enable us to push the argument through, however, we were forced to diverge from
Coleman's choice
of twist. We define the overconvergence of $F$, not in terms of
$F/E_{\kappa}$, but rather in terms of $F/h^s$ for some explicit 
modular form $h$. The benefit of our choice of $h$ is that it
is nicely compatible with the explicit formulae developed in~\cite{buzzcal},
and hence we may prove all our convergence results by hand in this case.
Our proof that eigenforms of finite slope overconverge ``as far as possible''
is essentially standard. The main contribution of this paper is to analyse the
overconvergence (or lack thereof)
of eigenforms in the kernel of the $U$ operator in this case.

One disadvantage of our approach is that the power series defining $h^s$
only converges
for $s$ sufficiently small and hence our arguments only deal with forms
whose weights lie in a certain disc at the centre of
weight space. However, recently in~\cite{buzzkill}, the 2-adic
level~1 eigencurve was shown to be a disjoint union of copies
of weight space near the boundary of weight space, and hence is automatically
proper here.

\section{Definitions}

Let~$\Delta(\tau)=q \prod_{n=1}^{\infty} (1 - q^n)^{24} = 
q - 24q^2 + 252q^3 - 1472q^4 +\ldots$ denote the classical
level~1 weight~12 modular form (where $q = e^{2 \pi i \tau}$). Set
$$f=\Delta(2 \tau)/\Delta(\tau)=q + 24q^2 + 300q^3 + 2624q^4 +\ldots,$$
a uniformizer for~$X_0(2)$, and
\begin{align*}
h=\Delta(\tau)^2/\Delta(2 \tau)&=\prod_{n\geq1}\left(\frac{1-q^n}{1+q^n}\right)^{24} = 1 - 48q + 1104q^2 - 16192q^3 +\ldots\\
\end{align*}
a modular form of level~2 and weight~12. Note that
the divisor of~$h$
is~$3(0)$, where~$(0)$ denotes the zero cusp on~$X_0(2)$, and hence that
$$h^{1/3} = \prod_{n\geq1}\left(\frac{1 - q^n}{1 + q^n}\right)^8$$
is a classical modular form of weight $4$ and level~2.

We briefly review the theory of overconvergent
$p$-adic modular forms, and make it completely explicit in the setting
we are interested in, namely $p=2$ and tame level~1. Let $\C_2$ denote
the completion of an algebraic closure of $\Q_2$. Normalise the norm
on $\C_2$ such that $|2|=1/2$, and normalise the valuation
$v:\C_2^\times\to\Q$ so that $v(2)=1$. Choose a group-theoretic
splitting of $v$ sending $1$ to $2$,
and let the resulting homomorphism $\Q\to\C_2^\times$
be denoted $t\mapsto 2^t$. Define $v(0)=+\infty$.
Let $\Ot_2$ denote the elements of $\C_2$ with non-negative valuation.

If~$r\in\Q$ with $0<r<2/3$ (note that $2/3=p/(p+1)$ if~$p=2$)
then there is a rigid space $X_0(1)_{\geq 2^{-r}}$ over~$\C_2$ such that
functions on this space are $r$-overconvergent 2-adic modular functions.
Let $X[r]$ denote the rigid space $X_0(1)_{\geq 2^{-r}}$.
By Proposition~1 of the appendix to~\cite{buzzcal},
we see that $X[r]$ is simply the closed
subdisc of the $j$-line defined by~$|j|\geq 2^{-12r}$.
We will also need to use (in Lemma~\ref{integercase}) the rigid
space $X[2/3]$, which we define as the closed subdisc of the
$j$-line defined by $|j|\geq 2^{-8}$.
The parameter~$q$ can be viewed as a rigid function defined in a neighbourhood
of~$\infty$ on $X[r]$, and hence any rigid
function on $X[r]$ can be written as a power series
in~$q$; this is the $q$-expansion of the form in this rigid analytic
setting. Moreover, it is well-known that the classical level~2
form~$f$ descends to a function on $X[r]$ (for any~$r<2/3$), with the
same $q$-expansion as that given above.

For $0<r<2/3,$ define $M_0[r]$
to be the space of rigid functions on $X[r]$, equipped with
its supremum norm. Then $M_0[r]$ is a Banach space over~$\C_2$ --- it is
the space of $r$-overconvergent modular forms of weight~0.
An easy calculation using the remarks after Proposition~1 of the appendix
to~\cite{buzzcal}
shows that the set
$\{1,2^{12r} \newf,2^{24r} \newf^2,\ldots,(2^{12r} \newf)^n,\ldots\}$
is an orthonormal Banach basis for~$M_0[r]$, and we endow $M_0[r]$
once and for all with this basis.

We define~$\W$ to be the open disc of centre~1 and radius~1 in the rigid
affine line over~$\C_2$. If $w\in\W(\C_2)$ then there is a unique continuous
group homomorphism $\kappa:\Z_2^\times\to\C_2^\times$ such that
$\kappa(-1)=1$ and $\kappa(5)=w$; moreover this establishes
a bijection between~$\W(\C_2)$ and the set of even 2-adic weights,
that is, continuous group homomorphisms
$\kappa:\Z_2^\times\to\C_2^\times$ such that $\kappa(-1)=1$. Note
that if~$k$ is an even integer then the map $x\mapsto x^k$ is
such a homomorphism, and we refer to this weight as weight~$k$.
Let $\tau:\Z_2^\times\to\C_2^\times$ denote the character
with kernel equal to $1+4\Z_2$, and let $\langle\cdot\rangle$ denote
the character $x\mapsto x/\tau(x)$; this character corresponds
to $w=5\in\W(\C_2)$. If $t\in\C_2$ with $|t|<2$ then
we may define $5^t:=\exp\bigl(t\log(5)\bigr)\in\W(\C_2)$ and we let
$\langle\cdot\rangle^t$ denote the homomorphism $\Z_2^\times\to\C_2^\times$
corresponding to this point of weight space. One
checks easily that the points of weight space corresponding to
characters of this form are $\{w\in\W(\C_2)\,:\,|w-1|<1/2\}$.

We now explain the definitions of overconvergent modular forms of
general weight that we shall use in this paper.
Recall $h=\prod_{n\geq1}(1-q^n)^{24}/(1+q^n)^{24}$.
Define $h^{1/8}$ to be the formal $q$-expansion
$\prod_{n\geq1}((1-q^n)^3/(1+q^n)^3$.
Now
$$(1-q^n)/(1+q^n)=1-2q^n+2q^{2n}-\ldots\in 1+2q\Z[[q]]$$
and hence $h^{1/8}\in 1+2q\Z[[q]]$.
Write $h^{1/8}=1+2qg$ with $g\in\Z[[q]]$. If $S$ is a formal variable
then we define $h^S\in 1+16qS\Z_2[[8S,q]]$ to be the formal
binomial expansion of $(1+2qg)^{8S}$. If $s\in\C_2$ with $|s|<8$
then we define $h^s$ to be the specialisation in $1+2q\Ot_2[[q]]$
of $h^S$ at $S=s$. In fact for the main part of this paper
we shall only be concerned with $h^s$ when $|s|<4$.

If $s\in\C_2$ with $|s|<8$, then define $\mu(s):=\min\{v(s),0\}$,
so $-3<\mu(s)\leq 0$. Define $\X$ to be the pairs $(\kappa,r)$ (where
$\kappa:\Z_2^\times\to\C_2^\times$ and $r\in\Q$) such that there
exists $s\in\C_2$ with $|s|<8$ satisfying

\begin{itemize}
\item $\kappa=\langle\cdot\rangle^{-12s}$, and
\item $0<r<1/2+\mu(s)/6$.
\end{itemize}

Note that the second inequality implies $r<1/2$, and conversely
if $|s|\leq1$ and $0<r<1/2$ then $(\langle\cdot\rangle^{-12s},r)\in\X$.

For $(\kappa,r)\in\X$,
and only for these $(\kappa,r)$,
we define the space $M_\kappa[r]$ of $r$-overconvergent forms
of weight $\kappa$ thus. Write $\kappa=\langle\cdot\rangle^{-12s}$
and define $M_\kappa[r]$ to be the
vector space of formal $q$-expansions $F\in\C_2[[q]]$ such that $Fh^s$
is the $q$-expansion of an element of $M_0[r]$.
We give $M_\kappa[r]$ the Banach space structure such that
multiplication by $h^s$ induces an isomorphism of Banach
spaces $M_\kappa[r]\to M_0[r]$, and we endow $M_\kappa[r]$ once
and for all with the orthonormal basis
$\{h^{-s},h^{-s}(2^{12r}f),h^{-s}(2^{12r}f)^2,\ldots\}$.

\begin{remark}\label{21} We do not consider the question here as to whether,
for all $(\kappa,r)\in\X$, the space
$M_\kappa[r]$ is equal to the space of $r$-overconvergent modular forms
of weight $\kappa$ as defined by Coleman (who uses the weight $\kappa$
Eisenstein series $E_\kappa$ to pass from weight $\kappa$ to weight~0).
One could use the methods of proof of \S5 of~\cite{buzzkill} to verify this;
the issue is verifying whether $E_\kappa h^s$ is $r$-overconvergent
and has no zeroes on $X[r]$. However, we do not need this result --- we
shall prove all the compactness results
for the $U$ operator that we need by explicit matrix computations,
rather than invoking Coleman's results. 
Note however that our spaces clearly coincide with Coleman's
if $\kappa=0$, as the two definitions coincide in this case.
Note also that for $r>0$ sufficiently small (depending on
$\kappa=\langle\cdot\rangle^{-12s}$ with $|s|<8$),
the definitions do coincide, because if $E_1:=1+4q+4q^2+\cdots$
denotes the weight~1 level~4 Eisenstein series, then
$h/E_1^{12}=1-96q+\cdots$ is overconvergent of weight~0, has
no zeroes on $X[r]$ for $r<1/3$, and has $q$-expansion congruent
to~1 mod~32. Hence for $r>0$ sufficiently small, the supremum
norm of $(h/E_1^{12})-1$ on $X[r]$ is $t$ with $t<1/2$ and $|s|t<1/2$, and this
is enough to ensure that the power series $(h/E_1^{12})^s$
is the $q$-expansion of a function on $X[r]$ with supremum norm at most~1.
Hence instead of using powers of $h$ to pass between weight $\kappa$
and weight~0, we could use powers of $E_1$. Finally, Corollary~B4.5.2
of~\cite{Coleman} shows that if $\kappa=\langle\cdot\rangle^{-12s}$
then there exists $r>0$ such that $E_1^{-12s}/E_\kappa$
is $r$-overconvergent, which suffices.
\end{remark}

Recall that 
if $X$ and $Y$ are Banach spaces over a complete field $K$ with orthonormal
bases $\{e_0,e_1,e_2,\ldots\}$ and $\{f_0,f_1,f_2,\ldots\}$,
then by the \emph{matrix} of a continuous linear map $\alpha:X\to Y$ we mean
the collection $(a_{ij})_{i,j\geq0}$ of elements of $K$ such
that $\alpha(e_j)=\sum_{i\geq0}a_{ij}f_i$. One checks that
\begin{itemize}
\item $\sup_{i,j}|a_{ij}|<\infty$, and
\item for all $j$ we have $\lim_{i\to\infty}|a_{ij}|=0$,
\end{itemize}
and conversely that given
any collection $(a_{ij})_{i,j\geq0}$ of elements of $K$ having
these two properties, there is a unique continuous linear map $\alpha:X\to Y$
having matrix $(a_{ij})_{i,j\geq0}$ (see Proposition~3 of~\cite{serre}
and the remarks following it for a proof).
When we speak
of ``the matrix'' associated to a continuous linear map between two
spaces of overconvergent modular forms, we will
mean the matrix associated to the map using the bases that
we fixed earlier.

If $R$ is a ring then we may define maps~$U$, $V$ and~$W$ on the ring
$R[[q]]$ by
\begin{align*}
U\left(\sum a_nq^n\right)&=\sum a_{2n}q^n,\\
V\left(\sum a_nq^n\right)&=\sum a_nq^{2n},\\
\mbox{and}\\
W\left(\sum a_nq^n\right)&=\sum (-1)^na_n q^n.
\end{align*}
Recall that $U(V(G)F)=GU(F)$ for $F,G$ formal power series in $q$,
and that $V:R[[q]]\to R[[q]]$ is a ring homomorphism.
The operator $W$ is not standard (or at least, our notation
for it is not standard),
but is also a ring homomorphism (it sends $f(q)$ to $f(-q)$)
and one also checks easily that $W=2VU-\Id$.
We shall show later on that there are continuous linear maps between various
spaces of overconvergent modular forms which correspond to $U$
and $W$, and will write down explicit formulae
for the matrices associated to these linear maps.

\section{The $U$ operator on overconvergent modular forms}

Our goal in this section is to make precise the statement in the
introduction that finite slope $U$-eigenforms overconverge a long
way. Fix $r\in\Q$ with $0<r<1/2$. We will show
that if $(\kappa,r)\in\X$ then the $U$-operator (defined on $q$-expansions)
induces a continuous linear map $M_\kappa[r]\to M_\kappa[r]$, and
we will compute the matrix of this linear map (with respect to our chosen basis
of $M_\kappa[r]$). We will deduce that if $0<\rho<r$ and $F$ is
$\rho$-overconvergent with $UF=\lambda F\not=0$ then $F$ is
$r$-overconvergent. These results are essentially standard but
we shall re-prove them, for two reasons: firstly to show that the
arguments still go through with our choice of twist, and secondly
to introduce a technique for computing matrices of Hecke operators
in arbitrary weight that we shall use when analysing the $W$ operator later.

It is well-known that the $U$-operator
induces a continuous linear map $U:M_0[r]\to M_0[r]$, and its associated
matrix was computed in~\cite{buzzcal}.
Now choose $m \in \Z_{\geq0}$, and set $k = -12m$.
One checks that $(k,r)\in\X$.
If $\phi\in M_0[r]$ then
$$h^m U \left(h^{-m} \phi\right) = 
h^m U\left(\Delta(2 \tau)^{-m} f^{2m}\phi\right) =
h^m \Delta(\tau)^{-m} U\left( f^{2m}\phi\right)
= f^{-m} U \left(f^{2m}\phi\right).$$
A simple analysis of the $q$-expansion of $f^{-m}U(f^{2m}\phi)$ shows
that it has no pole at the cusp of $X[r]$ and hence
$f^{-m}U(f^{2m}\phi)\in M_0[r]$. We deduce
that $U$ induces a continuous map $M_k[r]\to M_k[r]$,
and moreover that the matrix of this map (with respect to the basis fixed
earlier) equals
the matrix of the operator $U_k:=f^{-m} U f^{2m}$ acting on $M_0[r]$.
We now compute this matrix.

\begin{lemma}\label{s_is_m} For $m\in\Z_{\geq0}$ and $k=-12m$ as above,
and $j\in\Z_{\geq0}$, we have
$$U_k\left((2^{12r}f)^j\right)=\sum_{i=0}^{\infty}u_{ij}(m)(2^{12r}f)^i,$$
where $u_{ij}(m)$ is defined as follows: 
we have $u_{00}(0)=1$, $u_{ij}(m)=0$ if $2i-j<0$ or $2j-i+3m<0$, and
$$u_{ij}(m) = 
\frac{3(i+j + 3m - 1)! (j + 2m) 2^{8i - 4j + 12 r (j-i)}}
{2 (2i-j)! (2j-i + 3m)!}$$
if $2i-j\geq0$, $2j-i+3m\geq0$, and $i$, $j$, $m$ are not all zero.
\end{lemma}
\begin{proof}
The case $m=0$ of the lemma is Lemma~2 of~\cite{buzzcal},
and the general case follows easily from the fact that $U_k=f^{-m}U f^{2m}$.
Note that in fact all the sums in question are finite,
as $u_{ij}(m)=0$ for $i>2j+3m$.
\end{proof}

Now for $i,j\in\Z_{\geq0}$ define a polynomial $u_{ij}(S)\in\C_2[S]$ by
$u_{ij}(S)=0$ if $2i<j$, $u_{ij}(S)=2^{12 i r}$ if $2i=j$,
and 
$$u_{ij}(S) = \frac{3 \cdot 2^{12 r(j-i)} (j + 2S) 2^{8i-4j}}{2(2i-j)!} 
\prod_{\lambda=1}^{2i-j-1} (2j - i + \lambda + 3S)$$
if $2i > j$. One checks easily that evaluating $u_{ij}(S)$ at $S=m$
for $m\in\Z_{\geq0}$ gives $u_{ij}(m)$, so there is no
ambiguity in notation. Our goal now is to prove
that for all $s\in\C_2$ such that $|s|<8$ and
$(\langle\cdot\rangle^{-12s},r)\in\X$, the matrix
$(u_{ij}(s))_{i,j\geq0}$ is the matrix of the $U$-operator acting
on $M_\kappa[r]$ for $\kappa=\langle\cdot\rangle^{-12s}$ (with
respect to the basis of $M_\kappa[r]$ that we fixed earlier).

Say $s\in\C_2$ with $|s|<8$, define $\kappa=\langle\cdot\rangle^{-12s}$,
set $\mu=\min\{v(s),0\}$, and say $0<r<1/2+\mu/6$.
Then $(\kappa,r)\in\X$. 
Note that $v(as+b)\geq\mu$ for any $a,b\in\Z$, and
$3+\mu-6r>0$.
\begin{lemma} 

(a) One has
$v(u_{ij}(s)) \ge (3 + \mu - 6 r)(2i-j) + 6rj$.

(b) There is a continuous linear map $U(s):M_0[r]\to M_0[r]$ with
matrix $u_{ij}(s)$. Equivalently, there is a continuous linear map
$U(s):M_0[r]\to M_0[r]$ such that
$$U(s)\left((2^{12r}f)^j\right)=\sum_{i=0}^\infty u_{ij}(s)(2^{12r}f)^i.$$

\end{lemma}
\begin{proof}
(a) This is a trivial consequence of our explicit formula for $u_{ij}(s)$,
the remark about $v(as+b)$ above, and the fact that $v(m!) \le m - 1$
if $m\geq1$ (see Lemma~\ref{lemma:factorials}).

(b) Recall that $u_{ij}(s)=0$ if $2i<j$. Hence by (a) we see
that $|u_{ij}(s)|\leq 1$ for all $i,j$. It remains to check
that for all $j$ we have $\lim_{i\to\infty}v(u_{ij}(s))=+\infty$
which is also clear from (a).
\end{proof}
Note that $U(s)=U_{-12s}$ if $s=m\in\Z_{\geq0}$.

In fact the same argument gives slightly more. 
Choose $\epsilon\in\Q$ with $0<\epsilon<\min\{r,1/2+\mu/6-r\}$.
Then $(\kappa,r+\epsilon)\in\X$.
\begin{theorem} 
The endomorphism $U(s)$ of $M_0[r]$
is the composite of a continuous
map $M_0[r] \rightarrow M_0[r + \epsilon]$ and the restriction
$M_0[r + \epsilon] \rightarrow M_0[r]$.
\end{theorem}
\begin{proof} Define $w_{ij}(s)=u_{ij}(s)/2^{12\epsilon i}$.
By the previous lemma we have
$$v(w_{ij}(s))\geq(2i-j)(3+\mu-6r-6\epsilon)+6j(r-\epsilon)$$
and $w_{ij}(s)=0$ if $j>2i$. In particular $v(w_{ij}(s))\geq0$
for all $i,j$, and moreover
for all $j$ we have $\lim_{i\to\infty}w_{ij}(s)=0$.
The continuous linear map $M_0[r]\to M_0[r+\epsilon]$ with matrix
$(w_{ij}(s))_{i,j\geq0}$ will hence do the job.
\end{proof}
As usual say $|s|<8$, $\kappa=\langle\cdot\rangle^{-12s}$ and
$(\kappa,r)\in\X$.
\begin{cor}\label{cor:one}
The map $U(s):M_0[r]\to M_0[r]$ is compact and its characteristic
power series is independent of $r$ with $0<r<1/2+\mu/6$. Furthermore
if $0<\rho<r$ then any non-zero $U(s)$-eigenform with non-zero eigenvalue
on $M_0[\rho]$
extends to an element of $M_0[r]$.
\end{cor}
\begin{proof} This follows via standard arguments from the theorem;
see for example Proposition~4.3.2 of~\cite{Eigencurve}, although
the argument dates back much further.
\end{proof}

Keep the notation: $|s|<8$, $\kappa=\langle\cdot\rangle^{-12s}$,
$\mu=\min\{v(s),0\}$ and $0<r<1/2+\mu/6$, so $(\kappa,r)\in\X$.
We now twist $U(s)$ back to weight $\kappa$ and show that
the resulting compact operator is the $U$-operator (defined in
the usual way on power series).

\begin{prop}
\label{prop:35}
The compact endomorphism of $M_\kappa[r]$ defined
by $\phi\mapsto h^{-s}U(s)(h^s\phi)$ is the $U$-operator,
i.e., sends $\sum a_nq^n$ to $\sum a_{2n}q^n$.
\end{prop}
\begin{proof} It suffices to check the proposition for
$\phi=h^{-s}(2^{12r}f)^j$
for all $j\in\Z_{\geq0}$, as the result then follows by linearity.
If $S$ is a formal variable then recall that we may think of $h^S$
as an element of $1+16qS\Ot_2[[8S,q]])$
and in particular as an invertible element of $\Ot_2[[8S,q]]$.
Write $h^{-S}$ for its inverse.
We may think of $(h^S)U(h^{-S}(2^{12r}f)^j)$ as an element
of $\Ot_2[[8S,q]]$ (though not yet as an element of $M_0[r]$). Write
$$(h^S)U(h^{-S}(2^{12r}f)^j)=\sum_{i\geq0}\tilde{u}_{ij}(S)(2^{12r}f)^i$$
with $\tilde{u}_{ij}(S)\in\Ot_2[[8S]]\otimes\C_2$ (this is clearly possible
as $f=q+\ldots$). The proposition is just the statement that
the power series $\tilde{u}_{ij}(S)$ equals the polynomial $u_{ij}(S)$.
Now there exists some integer $N>>0$ such
that both $2^Nu_{ij}(S)$ and $2^N\tilde{u}_{ij}(S)$
lie in $\Ot_2[[8S]]$ (as $u_{ij}(S)$
is a polynomial).
Furthermore, Lemma~\ref{s_is_m} shows that $u_{ij}(m)=\tilde{u}_{ij}(m)$
for all $m\in\Z_{\geq0}$ and hence $2^N(u_{ij}(S)-\tilde{u}_{ij}(S))$
is an element of $\Ot_2[[8S]]$ with infinitely many zeroes in the
disc $|8s|<1$, so it
is identically zero by the Weierstrass approximation theorem.
\end{proof}

\begin{cor} If $(\kappa,r)\in\X$ and $\kappa=\langle\cdot\rangle^{-12s}$
then $U$ is a compact operator
on $M_\kappa[r]$ and its characteristic power series
coincides with the characteristic power series of $U(s)$ on $M_0[r]$.
Furthermore $F\in M_\kappa[r]$ is an eigenvector for $U$ iff
$Fh^s\in M_0[r]$ is an eigenvector for $U(s)$.
\end{cor}

\begin{proof} Clear.
\end{proof}

The utility of these results is that they allow us to
measure the overconvergence of a 
finite slope form $F$ of transcendental weight by instead
considering
the associated form $Fh^s$ in weight~$0$. This will be particularly
useful to us later on in the case when $F$ is in the kernel of~$U$.
We record explicitly what we have proved.
By an overconvergent modular
form of weight $\kappa$ we mean an element of $\bigcup_r M_\kappa[r]$,
where $r$ runs through the $r\in\Q$ for which $(\kappa,r)\in\X$.

\begin{cor}
\label{cor:converges} If $(\kappa,r)\in\X$ and $f$ is an overconvergent modular
form of weight $\kappa$ which is an eigenform for $U$ with non-zero
eigenvalue, then $f$ extends to an element of $M_\kappa[r]$.
\end{cor}
\begin{proof} This follows from \ref{cor:one} and \ref{prop:35}.
\end{proof}

In fact we will need a similar result for families of modular forms,
but our methods generalise to this case. We explicitly state
what we need.
\begin{cor}
\label{cor:convergesfam}
Let $A\subseteq\W$ be an affinoid subdomain, say $0<\rho<r<1/2$,
and assume that for all $\kappa\in A(\C_2)$ we have $(\kappa,r)\in\X$.
Let $F\in\Ot(A)[[q]]$ be an analytic family of $\rho$-overconvergent
modular forms, such that $UF=\lambda F$ for some $\lambda\in\Ot(A)^\times$.
Then $F$ is $r$-overconvergent.
\end{cor}
\qed

\section{The $W$ operator on overconvergent modular forms}

We need to perform a similar analysis to the previous section
with the operator $W$. Because $W=2VU-\Id$
we know that $W$ induces a continuous linear map
$V:M_0[r]\to M_0[r]$ for $r<1/3$ (for $r$ in this range,
$U$ doubles and then $V$ halves the radius of convergence).
Our goal in this section is to show that, at least for
$\kappa=\langle\cdot\rangle^{-12s}$ with $|s|<8$,  there
is an operator on weight $\kappa$ overconvergent modular forms
which also acts on $q$-expansions in this manner, and to compute
its matrix.

We proceed as in the previous section by firstly introducing
a twist of $W$. If $m\in\Z_{\geq0}$, if $k=-12m$
and if $\phi\in M_0[r]$ then the fact that $h(q)/h(-q)=(f(-q)/f(q))^2$
implies
$$h^mW(h^{-m}\phi)=f^{-2m}W(f^{2m}\phi)$$
and so we define the operator $W_k$ on $M_0[r]$, $r<1/3$,
by $W_k:=f^{-2m}W f^{2m}:M_0[r]\to M_0[r]$.

Set $g=Wf$, so $g(q)=f(-q)=-q + 24q^2 - 300q^3 +\ldots$. Because
$g=2VUf-f=48Vf+4096(Vf)^2-f$, we see that the $g$ can be regarded
as a meromorphic function on $X_0(4)$ of degree at most~4. Similarly
$f$ may be regarded as a function on $X_0(4)$ of degree~2.
Now the meromorphic function
$$(1+48f-8192f^2g)^2-(1+16f)^2(1+64f)$$
on $X_0(4)$ has degree at most~16 but the first 1000 terms
of its $q$-expansion can be checked to be zero on a computer,
and hence this function is identically zero.
We deduce the identity
$$g=\frac{1+48f-(1+16f)\sqrt{1+64f}}{8192f^2},$$
where the square root is the one of the form $1+32f+\ldots$,
and one verifies using the binomial theorem that
$g=\sum_{i\geq1}c_if^i$ with
\begin{align*}
c_i&:=(-1)^i2^{4i-4}\left(\frac{(2i+2)!}{(i+1)!(i+2)!}-\frac{(2i)!}{i!(i+1)!}\right)\\
&=(-1)^i2^{4(i-1)}\frac{3(2i)!}{(i-1)!(i+2)!}
\end{align*}
The other ingredient we need to compute the matrix of $W_k$
is a combinatorial lemma.
\begin{lemma} If $j\geq1$ and $i\geq j+1$ are integers
then
$$\sum_{a=j}^{i-1}\frac{3(2a+j-1)!j(2i-2a)!}{(a-j)!(a+2j)!(i-a-1)!(i-a+2)!}
=\frac{(2i+j)!(j+1)}{(i-j-1)!(i+2j+2)!}.$$
\end{lemma}
\begin{proof} Set $k=i-1-a$ and $n=i-1-j$ and then eliminate the variables
$i$ and $a$; the lemma then takes the form
$$\sum_{k=0}^nF(j,n,k)=G(j,n)$$
and, for fixed $n$ and $k$, both $F(j,n,k)$ and $G(j,n)$ are rational
functions of $j$. The lemma is now easily proved using Zeilberger's
algorithm (regarding $j$ as a free variable), which proves that the
left hand side of the equation satisfies an explicit (rather cumbersome)
recurrence relation of degree 1; however it is easily checked that the
right hand side is a solution to this recurrence relation, and this
argument reduces the proof of the lemma to the case $n=0$, where it
is easily checked by hand.
\end{proof}

We now compute the matrix of $W_k$ on $M_0[r]$ for $r<1/3$ and
$k=-12m$, $m\in\Z_{\geq0}$.
\begin{lemma} For $j\geq0$ we have
$$W_k\bigl(2^{12r}f)^j\bigr)=\sum_{i=0}^\infty \eta_{ij}(m)(2^{12r}f)^i,$$
where $\eta_{ij}(m)$ is defined as follows: we have $\eta_{ij}(m)=0$
if $i<j$, $\eta_{ii}(m)=(-1)^i$, and for $i>j$ we define
$$\eta_{ij}(m)=\frac{(2i + j-1 + 6m)! 3(j+2m) \cdot 2^{(4 -12 r)(i-j)} (-1)^i}
{(i-j)! (i+2j+6m)!}.$$
\end{lemma}
\begin{proof} We firstly deal with the case $m=0$, by induction on $j$.
The case $j=0$
is easily checked as $\eta_{i0}(0)=0$ for $i>0$, and the case $j=1$
follows from the fact that $c_i2^{12r(1-i)}=\eta_{i1}(0)$ for $i\geq1$, as is
easily verified. For $j\geq1$ we have $W(f^{j+1})=f(-q)^{j+1}
=g\cdot W(f^j)=(\sum_{t\geq 1}c_tf^t)W(f^j),$
and so to finish the $m=0$ case it suffices to verify that for $j\geq1$
and $i\geq j+1$ we have
$\eta_{i\, j+1}(0)=2^{12r}\sum_{a=0}^{i-1}c_{i-a}2^{-12r(i-a)}\eta_{aj}(0),$
which quickly reduces to the combinatorial lemma above.

Finally we note that because $\eta_{i+2m\,j+2m}(0)=\eta_{ij}(m)$,
the general case follows easily from the case $m=0$ and the
fact that $W_k=f^{-2m}Wf^{2m}$.
\end{proof}
As before, we now define polynomials $\eta_{ij}(S)$ by $\eta_{ij}(S)=0$
if $i<j$, $\eta_{ii}(S)=(-1)^i$, and
$$\eta_{ij}(S) = \frac{3 (j + 2 S) 2^{(4 - 12r)(i-j)} (-1)^i}{(i-j)!}
\prod_{\lambda = 1}^{i-j-1}
(i + 2j + \lambda + 6S)$$
for $i>j$. We observe that $\eta_{ij}(S)$ specialises to $\eta_{ij}(m)$
when $S=m\in\Z_{\geq0}$. 

Now if $|s|<8$ and $\kappa=\langle\cdot\rangle^{-12s}$, and we
set $\lambda=\min\{v(2s),0\}>-2$, then we check easily
that $v(\eta_{ij}(s))\geq(3-12r+\lambda)(i-j)+1$,
so for $12r<3+\lambda$ we see that $(\eta_{ij}(s))_{i,j\geq0}$ is the matrix
of a continuous endomorphism $W(s)$
of $M_0[r]$. Moreover, arguments analogous to those of the previous
section show that if furthermore $(\kappa,r)\in\X$ (so $M_\kappa[r]$
is defined), then the
endomorphism of $M_\kappa[r]$ defined by sending
$\phi$ to $h^{-s}W(s)(h^s\phi)$ equals the $W$ operator as defined
on $q$-expansions. Note that if $|s|\leq 4$ then $12r<3+\lambda$
implies $(\kappa,r)\in\X$.

%{\bf SUMMARY} If $0<r<2/3$ and $|s|<8$ and $\kappa=\langle\cdot\rangle^{-12s}$
%then we have defined an operator $U(s)$ if $0<r<1/2+\mu/6$ (where
%$\mu=\min\{v(s),0\}$) and an operator $W(s)$ if $0<r<1/4+\lambda/12$
%where $\lambda=\min\{v(2s),0\}=\min\{\mu+1,0\}$. In fact if $d(2s,\Z_2)$ is
%less than 1 then we might be able to do better with $W(s)$ and if
%$d(s,\Z_2)$ is less than 1 then we might be able to do better with $U(s)$
%(I mean, define it for a few more $r$).

\section{Strategy of the proof.}

We have proved in Corollary~\ref{cor:converges} that overconvergent
modular forms $f$ such that $Uf=af$ with $a\not=0$ overconverge
``a long way''. Using the $W$-operator introduced in the previous
section we will now prove that overconvergent modular forms $f=q+\ldots$
such that $Uf=0$ cannot overconverge as far. We introduce a definition
and then record the precise statement.

\begin{df}\label{51}
 If $x\in\C_2$ then
set $\beta=\beta(x)=\sup\{v(x-n):n\in\Z_2\}$, allowing $\beta=+\infty$ if
$x\in\Z_2$,
and define $\nu=\nu(x)$ as follows:
$\nu=\beta$ if $\beta\leq0$,
$\nu=\beta/2$ if $0\leq\beta\leq1$,
and in general
$$\nu=\sum_{k=1}^n1/2^k+(\beta-n)/2^{n+1}$$
if $n\leq\beta\leq n+1$. Finally define $\nu=1$ if $\beta=+\infty$.
\end{df}

The meaning of the following purely elementary lemma will become
apparent after the statement of Theorem~\ref{52}.

\begin{lemma}\label{lemma:smalldisc} Say $s\in\C_2$ with $|s|<4$ and
furthermore assume $2s\not\in\Z_2^\times$.
Then for all $s'\in\C_2$
with $|s-s'|\leq 1$, we have
$0<\frac{3+\nu(2s)}{12}<\frac{1}{2}+\frac{\mu(s')}{6}.$
\end{lemma}
\begin{proof} We have $\nu(2s)>-1$ and so certainly $\frac{3+\nu(s)}{12}>0$.
The other inequality can be verified on a case-by-case basis.
We sketch the argument. 

If $|s|>2$ then $|s'|=|s|>2$
and $\nu(2s)-1=v(s)=v(s')=\mu(s')$; the
inequality now follows easily from the fact that $\mu(s')>-2$.

If $|s|\leq2$ but $2s\not\in\Z_2$ then $0<\beta(2s)<\infty$ and
$\nu(2s)<1$; now $|s'|\leq 2$ and hence $\mu(s')\geq-1$,
thus $\frac{3+\nu(2s)}{12}<\frac{1}{3}\leq\frac{1}{2}+\frac{\mu(s')}{6}$.

Finally if $2s\in\Z_2$ then we are assuming $2s\not\in\Z_2^\times$
and hence $s\in\Z_2$ so $|s|\leq1$ and hence $|s'|\leq 1$.
Hence $\mu(s')=0$ and we have
$\frac{3+\nu(2s)}{12}=\frac{1}{3}<\frac{1}{2}+\frac{\mu(s')}{6}$.
\end{proof}

Again say $|s|<4$ and $2s\notin\Z_2^\times$. Write
$\kappa=\langle\cdot\rangle^{-12s}$, and $\nu=\nu(2s)$.
Let $G=q+\ldots$ be an overconvergent form of weight $\kappa$
(by which we mean an element of $M_\kappa[\rho]$ for some $\rho\in\Q_{>0}$
sufficiently small). 
The theorem we prove in the next section (which is really the main
contribution of this paper) is

\begin{theorem}\label{52}
If $G=q+\ldots$ satisfies $UG=0$, then $F:=h^{s}G\in M_0[\rho]$
does not extend
to an element of $M_0[r]$ for $r=\frac{3+\nu}{12}$. Equivalently,
$G\not\in M_\kappa[r]$.
\end{theorem}

Note that by Lemma~\ref{lemma:smalldisc} we have $(\kappa,r)\in\X$
so the theorem makes sense.
Furthermore, by Corollary~\ref{cor:converges},
overconvergent eigenforms of the form $q+\ldots$ in the kernel
of $U$ overconverge less than finite slope overconvergent eigenforms.
Note also that if $2s\in\Z_2^\times$ then $\nu(2s)=1$ and
for $\kappa,r$ as above we have $(\kappa,r)\not\in\X$.
We deal with this minor annoyance in the last section of this paper.

\section{The Kernel of $U$}

In this section we prove Theorem~\ref{52}. We divide the argument
up into several cases depending on the value of $s$.
We suppose that $|s|<4$ and $2s\not\in\Z_2^\times$,
and we set $\kappa=\langle\cdot\rangle^{-12s}$.
Define $\nu=\nu(2s)$ as in the previous section, and set $r=\frac{3+\nu}{12}$.
For simplicity we drop the $s$ notation from $\eta_{ij}(s)$ and write
\begin{align*}
\eta_{ij}&=\frac{3(j+2s)2^{(4-12r)(i-j)}(-1)^i}{(i-j)!}
\prod_{t=1}^{i-j-1}(i+2j+t+6s)\\
&=\frac{3(j+2s)2^{(1-\nu)(i-j)}(-1)^i}{(i-j)!}
\prod_{t=1}^{i-j-1}(i+2j+t+6s).
\end{align*}

Say $G=q+\ldots$ as in Theorem~\ref{52} is $\rho$-overconvergent
for some $0<\rho<r$, so $F=h^sG\in M_0[\rho]$.
If we expand $F$ as
$$F=\sum_{j\geq 1}\tilde{a}_j(2^{12\rho}f)^j$$
then it follows that $\tilde{a}_1\not=0$. Recall
also that $\tilde{a}_j\to 0$ as $j\to\infty$.
On the other hand, $F = - W(s) F$,
and so
$$\tilde{a}_i = - \sum_{j=1}^{\infty} \tilde{a}_j \tilde{\eta}_{i,j},$$
where $\tilde{\eta}_{ij}$ denotes the matrix of $W(s)$ on $M_0[\rho]$
(so $\eta_{ij}=\tilde{\eta}_{ij}2^{12(r-\tilde{r})(j-i)}$).
We deduce from this that if we define
$a_i=\tilde{a}_i2^{12(\tilde{r}-r)i}$
then
$F=\sum_{j\geq1}a_j(2^{12r}f)^j$
and
$$a_i=-\sum_{j\geq 1}a_j\eta_{ij}.$$
Note in particular that the sum converges even if $W(s)$ does
not extend to a continuous endomorphism of $M_0[r]$ or if $F$
does not extend to an element of $M_0[r]$. In fact our goal is
to show that the $a_i$ do not tend to zero, and in particular
that $F$ does not extend to an element of $M_0[r]$.
\begin{lemma}  \label{lemma:utility}
Suppose $F$ is as above.
Suppose also that there exist constants $c_1$ and $c_3\in\R$,
an infinite set $I$ of positive integers, and
for each $i\in I$ constants $N(i)$ and $c_2(i)$ tending to infinity
as $i\to\infty$ and such that
\begin{itemize}
\item[(i)] $v(\eta_{i1}) \le c_1$, for all $i\in I$.
\item[(ii)] $v(\eta_{ij}) \ge c_2(i)$ for all  $i\in I$
and $2 \le j \le N(i)$.
\item[(iii)] $v(\eta_{ij}) \ge c_3$  for all $i\in I$ and $j\in\Z_{\geq0}$.
\end{itemize}
Then the $a_i$ do not tend to zero as $i\to\infty$,
and hence $F$ does not extend to a function on
$M_0[r]$.
\end{lemma}

\begin{proof} Assume $a_i\to 0$. Recall that we assume $a_1\not=0$.
By throwing away the first
few terms of $I$ if necessary, we may then assume that for all $i\in I$
we have
\begin{enumerate}
\item[(1)] $c_2(i)>v(a_1)+c_1-\min\{v(a_j):j\geq1\}$, and
\item[(2)] $\min\{v(a_j):j>N(i)\}>v(a_1)+c_1-c_3$.
\end{enumerate}
We now claim that for all $i\in I$
we have $v(a_1\eta_{i1})<v(a_j\eta_{ij})$ for all $j>1$.
The reason is that if $j\leq N(i)$ the inequality follows from
equation (1) above,
and if $j>N(i)$ it follows from (2).
Now from the equality
$$a_i = - \sum_{j=1}^{\infty} a_j \eta_{ij}$$
we deduce that $v(a_i)=v(a_1\eta_{i1})$ is bounded
for all $i\in I$, contradicting the fact that $a_i\to 0$.
\end{proof}

The rest of this section is devoted to establishing these inequalities
for suitable $I$ and $r$. We start with some preliminary lemmas.

\begin{lemma}
\label{lemma:factorials}
\begin{enumerate}
\item If $m\geq1$ then $v(m!)\leq m-1$, with equality
if and only if $m$ is a power of~2.
\item If $m\geq0$ then $v(m!)\geq (m-1)/2$, with equality if and only
if $m=1,3$.
\item If $n\geq0$ and $0\leq m<2^n$ then setting $t=2^n-m$
we have $m-v(m!)\geq n-(t/2)$.
\end{enumerate}
\end{lemma} 
\begin{proof}
1 and 2 follow easily from
$$v(m!)=\lfloor m/2\rfloor+\lfloor m/4\rfloor+\lfloor m/8\rfloor+\ldots.$$
For 3, we have $m!(m+1)(m+2)\ldots(2^n-1)(2^n)=(2^n)!$
and for $0<d<2^n$ we have $v(d)=v(2^n-d)$, so
$v((m+1)(m+2)\ldots(2^n-1))=v((t-1)!)\geq(t-2)/2$ by 2.
Finally $v((2^n)!)=2^n-1$ by 1.
Hence $v(m!)\leq2^n-1-n-(t-2)/2=2^n-n-(t/2)$ and so
$m-v(m!)\geq 2^n-t-(2^n-n-(t/2))=n-(t/2)$.
\end{proof}

\begin{lemma} \label{lemma:question}  Let $m\in\Z$ be arbitrary
and set $\beta=\beta(x)$ and $\nu=\nu(x)$ as in Definition~\ref{51}.
\begin{enumerate}
\item If $\beta\leq0$ then $v(x+n)=\nu$ for all $n\in\Z$, hence
the valuation of $\prod_{t=1}^N(x+m+t)$ is $N\nu$.
\item  If $0<\beta<\infty$ and if $N$ is a power of~2 with
$N\geq2^{\lceil\beta\rceil}$
then the valuation of 
$$\prod_{t=1}^N(x+m+t)$$
is exactly $N\nu$.
\item If $0<\beta<\infty$ and if $N\geq0$ is an
arbitrary integer then the valuation of  
$$\prod_{t=1}^N(x+m+t)$$
is $v$, where $|v-N\nu|<\beta$.
\item If $\beta=\infty$ and if $N\geq0$ is an arbitrary integer
then the valuation of $\prod_{t=1}^N(x+m+t)$ is at least $v(N!)$.
\end{enumerate}
\end{lemma}

\begin{proof} (1) is obvious and (2) is easy to check (note that $v(x+n)$
is periodic with period $2^{\lceil\beta\rceil}$). For part (3), 
say $n=\lfloor\beta\rfloor$.
Now about half
of the terms in this product are divisible by $2$, about a quarter
are divisible by $4$, and so on. More precisely, this means that 
the largest possible power of $2$ that can divide this product is
\begin{align*}
&\lceiling N/2 \rceiling + \lceiling N/4 \rceiling + \ldots + \lceiling N/2^n \rceiling + (\beta - n) \lceiling N/{2^{n+1}} \rceiling\\
<& (N/2+1)+(N/4+1)+\ldots+(N/2^n+1)+(\beta-n)(N/2^{n+1}+1)\\
=&N \nu + \beta.
\end{align*}

A similar argument shows that the lowest possible power of~2
dividing this product is strictly greater than $N\nu-\beta$.

For part (4), if $\beta=\infty$ then $x\in\Z_2$ and by a continuity
argument it suffices to prove the result for $x$ a large positive
integer, where it is immediate because the binomial coefficient
$\binom{x+m+N}{N}$ is an integer.

\end{proof} 

Now set $x=2s$ and let $\beta=\beta(2s),\nu=\nu(2s)$.
Note that if $\beta\leq0$ then $\mu=\beta-1$, and if $\beta\geq1$
then $\mu=0$.

Recall $\eta_{ij}=0$ if $i<j$, $\eta_{ii}=(-1)^i$,
and if $i>j$ we have 
$$\eta_{ij} = \frac{3 (j + 2 s) 2^{(1 - \nu)(i-j)} (-1)^i}{(i-j)!}
\prod_{t = 1}^{i-j-1}
(i + 2j + t + 6s).$$
In particular, for $i>j$ we have
$$(*)\ \ \ v(\eta_{i,j}) = (1-\nu)(i-j)-v((i-j)!) + v(j + 2s)
+ v \left( \prod_{t=1}^{i-j-1} (i + 2j + t + 6s) \right).$$
We shall continually refer to $(*)$ in what follows.
\begin{prop} Say $\beta\leq0$ (and hence $\nu=\beta$).
\label{prop:buzz1}
\begin{enumerate}
\item If $j\geq i$ then $v(\eta_{ij})\geq0$, and if $j<i$ then
$v(\eta_{ij})=i-j-v((i-j)!)\geq1$.
\item
 If $i=2^n+1$ then $v(\eta_{i1})=1$ and if $1<j<i$
then $v(\eta_{ij})\geq n-(j-1)/2$.
\end{enumerate}
\end{prop}

\begin{proof} 1 is immediate from $(*)$ and
Lemma~\ref{lemma:question}(1). Now 2 can be deduced from 1,
using part 1 of Lemma~\ref{lemma:factorials} for the first
part and part 3 of Lemma~\ref{lemma:factorials} for the second.
\end{proof}

\medskip

We now prove:
\begin{lemma} 
\label{lemma:two} Theorem~\ref{52} is true if $-1<\beta\leq0$ (i.e.,
if $2\leq|s|<4$).
\end{lemma}
Equivalently, if $2\leq|s|<4$ and $\kappa=\langle\cdot\rangle^{-12s}$,
and if $G=q+\ldots$ is a non-zero weight $\kappa$ overconvergent
form in $\ker(U)$, then $F=h^sG$ does not converge as far as
$M_0[1/4+\nu/12]$, where $\nu=\nu(2s)$ as above.

\begin{proof} This will be a direct application
of lemma~\ref{lemma:utility}.
We set $I=\{2^n+1:n\in\Z_{>0}\}$,
and if $i=2^n+1$ we define $c_2(i)=(n+1)/2$ and $N(i)=n$.
We set $c_1 = 1$ and $c_3 = 0$. Now assumptions (i) and (ii)
of Lemma~\ref{lemma:utility} follow from Proposition~\ref{prop:buzz1}(2),
and (iii) follows from Proposition~\ref{prop:buzz1}(1).
\end{proof}

\medskip

Let us now consider the case when $0<\beta<\infty$.

\begin{prop}
\label{57}
 Let  $0<\beta<\infty$. 
\begin{enumerate}
\item If $j<i$ then $v(\eta_{i,j})-\bigl((i-j)-v((i-j)!)-\nu\bigr)\in[-\beta,2\beta]$.
\item If $j<i$ then
$$v(\eta_{ij})\geq1-\beta-\nu.$$
If $i=2^n+1$ then
$$v(\eta_{i1})\leq2\beta-\nu+1$$
and if $1<j<i$ then
$$v(\eta_{ij})\geq n-(j+1)/2-\nu-\beta.$$
\end{enumerate}
\end{prop}
\begin{proof}
From the definition of $\beta$, the valuation of $j + 2s$
lies in $[0,\beta]$. The result then follows
from $(*)$ and lemma~\ref{lemma:question}, part 3.
Part 2 follows from part 1 and Lemma~\ref{lemma:factorials}, parts
(1) and (3), applied to $(i-j)!$.
\end{proof}

\begin{lemma} \label{lemma:three} Theorem~\ref{52}
is true if $0<\beta< \infty$, that is, if $|s|\leq 2$ and $2s\not\in\Z_2$.
\end{lemma}

\begin{proof}
Again this is an application of lemma~\ref{lemma:utility}.
Set $I=\{2^n+1:n\in\Z_{>0}\}$, 
$c_1=2 \beta - \nu+1$, $c_3 =\min\{0,1 - \beta-\nu\}$,
and if $i=2^n+1$ then set $N(i)=n$ and $c_2(i)=(n+1)/2 - \nu - \beta$.
Conditions (i)--(iii) of Lemma~\ref{lemma:utility}
hold by Proposition~\ref{57}(2).
\end{proof}

The only cases of Theorem~\ref{52} left to deal with are those
with $\beta=+\infty$, that is, $2s\in\Z_2$. Because the theorem
does not deal with the case $2s\in\Z_2^\times$ we may assume
from now on that $2s\in 2\Z_2$, so $s\in\Z_2$. We next deal with
the case $s\in\Z_2$ and $6s\not\in\N$, where $\N=\{1,2,3,\ldots\}$
is the positive integers.
In this case, we shall again use Lemma~\ref{lemma:utility} with $i$
of the form $i = 2^n +1$. However, it will turn out that only certain
(although infinitely many) $n$ will be suitable.

Since we assume $s\in\Z_2$ we have $\beta=+\infty$, so $\nu=1$ and hence
$$(**)\ \ \ \eta_{ij} = \frac{3 (j + 2 s)(-1)^i}{(i-j)!}
\prod_{t = 1}^{i-j-1}
(i + 2j + t + 6s) .$$

Let $u \in \Z_2$. Define functions $f_n(u)$ as follows:
$$f_n(u) = (2^n + u)(2^{n} + u + 1) \cdots (2^{n+1} - 1 + u) = 
\prod_{\tau=0}^{2^n - 1}(2^n + u + \tau).$$
\begin{lemma}\label{lemma:68}
For any $u \in \Z_2$ there exist infinitely many values of $n$
for which
$$v_n(f(u)) = v((2^n)!) \ \text{or} \ v((2^n)!) + 1.$$
\end{lemma}
\begin{proof} For each $n$, define an integer $0 < u_n \le 2^n$ by setting
$u \equiv u_n \mod 2^n$.
If $0\leq \tau\leq 2^n-1$ and $\tau \ne 2^n - u_n$, then
$$v(2^n + u + \tau) = v(u_n + \tau).$$
Since $\tau$ takes on every equivalence class modulo $2^n$,
It follows from the definition of $f_n$ that
$$v(f_n(u)) = v((2^n - 1)!) + v(2^{n+1} + u - u_n).$$
If $u \not \equiv u_n \mod 2^{n+1}$ then $v(2^{n+1} + u - u_n)
= v(2^n)$ and $v(f_n(u)) = v((2^n)!)$.
There are infinitely many $n$ satisfying this condition
unless $u \equiv u_n \mod 2^{n+1}$ for all sufficiently
large~$n$. Yet
this implies $u_n = u_{n+1}$ for all sufficiently large $n$,
and subsequently that $u = u_n$. In this case
we have $v(2^{n+1} + u - u_n) =
v(2^{n+1})$, and $v(f_n(u)) = v((2^n)!) + 1$.
\end{proof}

\begin{cor}\label{69}
There are infinitely many $n$ such that if $i=2^n+1$
then $v(\eta_{i1})\in\{0,1\}$.
\end{cor}
\begin{proof} Let $i=2^n+1$ and $j=1$, and assume $n\geq1$.
By $(**)$ we have
$$\eta_{i1} = \frac{3 (1 + 2 s)(-1)}{(2^n)!}
\prod_{t = 1}^{2^n-1}
(2^n + 3 + t + 6s) .$$
Let $u = 6s+4\in2\Z_2$
and set $\tau=t-1$. Then
$$\eta_{i1} = \frac{(1 - u)}{(2^n)!}
\prod_{\tau = 0}^{2^n - 2}
(2^n + u + \tau) 
=  \frac{f_n(u)}{(2^n)!} \cdot \frac{1 - u}{u - 1 + 2^{n+1}}$$
and the result follows from Lemma~\ref{lemma:68}
and the fact that $u\in2\Z_2$.
\end{proof}

Let us now turn to estimating $\eta_{ij}$ for general $i,j$.

\begin{lemma}
\label{combinatorics}
If $i,j\in\Z_{\geq0}$ then $v(\eta_{ij})\geq0$.
\end{lemma}
\begin{proof}By continuity, it suffices to verify the result
for $6s$ a large positive even integer. It is clear if $i\leq j$
so assume $i>j$. Now because the product of $N$
successive integers is divisible by $N!$ we see (putting one
extra term into the product) that both $x_1:=\frac{i+2j+6s}{3(j+2s)}\eta_{ij}$
and $x_2:=\frac{2i+j+6s}{3(j+2s)}\eta_{ij}$ are integers.
The result now follows as $\eta_{ij}=2x_1-x_2$.
\end{proof}

Set $I_0=\{i=2^n+1:v(\eta_{i1})\in\{0,1\}\}$. Then $I_0$ is infinite
by Corollary~\ref{69}. We will ultimately let $I$
be a subset of $I_0$. We must analyse $\eta_{ij}$ for $i\in I_0$ and
$1<j$ small.
Note that if $i=2^n+1$ and $j\geq2$, then
$$\frac{\eta_{i,j}}{\eta_{i,1}} =  2^n \cdot \frac{(j + 2s)}{(1 + 2s)}
\cdot \prod_{t=1}^{j-2} (i - j + t) 
\cdot \frac{ \prod_{t=1}^{j-1} (2i + t + 6s)}
{\prod_{t=1}^{2j-2} (i + 2 + t + 6s)}$$
Since $6s \notin -\N$, $3 + 6s + t \neq 0$. Thus
for any $N$ there exists $n_0$ depending on $N$ such
that for all $n\geq n_0$ we have
$v(i + 2 + 6s + t) = v(3 + 6s + t)$ for all $t\leq 2N-2$.
In particular, for
fixed $N$ and sufficiently large $n$ (with $i = 2^n + 1$),
$$v(\eta_{ij}) \ge n  - v\left( \prod_{t=0}^{2j-2} (3 + 6s + t)\right)
 + v(\eta_{i1}).$$
\begin{lemma}\label{611}
For any constants $c_2\in\R$ and $N\in\Z_{\geq1}$,
there exists $n_1=n_1(c_2,N)$ such that for all $n\geq n_1$
such that $i=2^n+1\in I_0$, we have
$v(\eta_{ij}) \ge c_2$ for $2 \le j \le N$.
\end{lemma}
\begin{proof}
Set $M=v(\prod_{t=0}^{2N-2} (3 + 6s + t))$
and choose $n_1$ such that $n_1-M\geq c_2$.
\end{proof}

\medskip

We may now prove:
\begin{lemma} 
\label{lemma:four} 
Theorem~\ref{52} is true if $s\in\Z_2$ and $6s\not\in-\N$.
\end{lemma}

\begin{proof}
We apply lemma~\ref{lemma:utility} as follows.
Set $c_1=1$ and $c_3=0$. We build $I$ as follows.
As $m$ runs through the positive integers, set $N=c_2=m$,
define $n_1$ as in Lemma~\ref{611}, choose $n\geq n_1$ such that
$i:=2^n+1\in I_0$ and such that $i$ is not yet in $I$; now add $i$
to $I$ and define $N(i)=c_2(i)=t$.
The conditions of lemma~\ref{lemma:utility}
are then satisfied.
\end{proof}

The final case in our proof of Theorem~\ref{52} is the
case $6s \in -2\N$, which corresponds to weight 
$k = -12s \in 4 \N$.
We shall not
use Lemma~\ref{lemma:utility} in this case, but give a direct argument.

Because our level structure is so small it is convenient to temporarily
augment it to get around representability issues.
So choose some odd prime $p$ and consider the compact modular curve~$Y$
over $\Q_2$ whose cuspidal points parametrise elliptic curves
with a subgroup of order~2 and a full level $p$ structure (note that this
curve is not in general connected). There is a sheaf~$\omega$
on~$Y$, and classical modular forms of weight $k$ and level~2
are, by definition, $\GL_2(\Z/p\Z)$-invariant global
sections of $\omega^{\otimes k}$ on $Y$.

For $0<r\leq 2/3$ let $Y[r]$ denote the pre-image of $X[r]$ via the
forgetful functor. Recall that there is a compact
operator $U$ on $H^0(Y[r],\omega^{\otimes k})$ for $r<2/3$ and $k\in\Z$.

\begin{lemma}
\label{integercase}
If $k\in\Z$ and $f\in H^0(Y[1/3],\omega^{\otimes k})$
is in the kernel of $U$, then $f=0$.
\end{lemma}

\begin{remark} The lemma is not special to $p=2$; the proof shows
that non-zero $p$-adic modular forms in the kernel of $U$ are never
$1/(p+1)$-overconvergent.
\end{remark}

\begin{proof} Say $f\in H^0(Y[1/3],\omega^{\otimes k})$ is
arbitrary. If $E$ is an elliptic curve over a finite extension
of $\Q_2$, equipped with with a subgroup $C$
of order $2$ and a full level $p$ structure $L$,
and such that the corresponding point $(E,C,L)\in Y$ is in $Y[1/3]$,
then one can regard $f(E,C,L)$ as an element
of $H^0(E,\Omega^1)^{\otimes k}$. Now define
$g\in H^0(Y[2/3],\omega^{\otimes k})$ by
$$g(E,L)=\sum_{D\not=C}(\pr)^*f(E/D,\overline{C},\overline{L}),$$
where the sum is over the subgroups~$D\not=C$ of $E$ of order~2,
$\pr$ denotes the projection $E\to E/D$, and a bar over a level
structure denotes its natural pushforward. An easy calculation
using Tate curves (see for example Proposition~5.1 of~\cite{wild})
shows that $g=pUf$, and hence if $Uf=0$ then $g=0$. In particular
if $E$ is an elliptic curve with no canonical subgroup and we fix
a full level $p$ structure $L$ on $E$, then
then $(E,C,L)\in Y[2/3]$ for all $C$, and $g(E,C,L)=0$
for all $C$ implies that $\sum_{D\not=C}(\pr)^*f(E/D,E[2]/D,\overline{L})=0$
for all $C$. Summing, one deduces that
$\sum_D(\pr)^*f(E/D,E[2]/D,\overline{L})=0$ and hence that
$f(E/D,E[2]/D,\overline{L})=0$ for all $D$ of order~2.
This implies that $f$ is identically zero on the ``boundary''
of $Y[1/3]$ and hence that $f$ is identically zero.
\end{proof}

We deduce

\begin{lemma} \label{lemma:five} Theorem~\ref{52} is true for $6s \in -2 \N$.
\end{lemma}
\begin{proof} If $G\in M_k[1/3]$ then $G=h^{k/12}F$ and, because
$k=-12s\in 4\N$, we know that $h^{k/12}$ is a classical modular form
of level~2 and hence an element of $H^0(Y[1/3],\omega^{\otimes(k/12)})$.
Thus the preceding lemma applies to $G$ and we conclude that $G=0$.
\end{proof}

Theorem~\ref{52} now follows from Lemmas~\ref{lemma:two}, \ref{lemma:three},
\ref{lemma:four} and \ref{lemma:five}.

\section{General facts about the 2-adic eigencurve.}

In this section we collect some standard results about $\E$,
including several for which we know no reference. For brevity
we have restricted our attention to the 2-adic level~1 eigencurve,
but much of what we say applies more generally (see Remark~\ref{76} for
more precise comments about what works in general and what doesn't.) We
remark that sometimes our proofs could be shortened slightly but
we have presented proofs that would generalise easily once one has
set up the required notation.

We firstly recall the definition of the eigencurve $\E$,
following Part~II of~\cite{eigenvarieties}.
If $Y=\Sp(R)\to\W$ is a map from an affinoid to $\W$
(for example a point of $\W$ or an admissible affinoid
open in $\W$) then for $0<r<2/3$ we define
$M_Y[r]$ to be the space of $r$-overconvergent modular
forms of weight~$Y$, that is, the $R$-module of formal power series
$F\in R[[q]]$ such that $F/E_Y$ is the $q$-expansion
of an element of $\Ot(X[r]\times Y)$, where
$E_Y$ denotes the pullback of the Eisenstein family in $\Ot(\W)[[q]]$
to $R[[q]]$.
This is probably not the ``correct'' definition if $r$ is
close to $2/3$ and the image of $Y$ contains points near the
boundary of weight space. On the other hand, it is shown
in section~7 of~\cite{eigenvarieties} that for $r$ sufficiently
close to zero the space $M_Y[r]$ is stable under all the Hecke
operators $T_n$ (this was not proved in~\cite{Eigencurve}, although
it was stated for $p>2$; however the missing ingredient is provided
by Lemma~7.1 of~\cite{eigenvarieties}). We assume henceforth
that $r$ is always sufficiently small for all the Hecke operators
to be defined. Then $U$ is a well-defined compact
endomorphism of $M_Y[r]$ and its characteristic power series
$P_{U,Y}(T)\in R[[T]]$
is independent of $r$. As noted in Remark~\ref{21},
if the image of $Y$ in weight space is contained
within the characters of the form $\langle\cdot\rangle^{-12s}$ with
$|s|<8$, and if $r$ is sufficiently small, then the
definition of $M_Y[r]$ above coincides with the one used in this paper.
We henceforth assume that $r$ is also sufficiently small to ensure
that this is the case.

Now if $\phi$ is a compact operator on $M_Y[r]$
(for example, $\phi$ could be the operator $U$, or $UT_\ell$ for some
odd prime $\ell$) then its characteristic power series $P_{\phi,Y}(T)$ is in
$R[[T]]$ and can
even be regarded as a function on $Y\times\A^1$. Let $Z_{\phi,Y}$ denote
the closed subspace of $Y\times\A^1$ cut out by $P_{\phi,Y}(T)$.
If we let $Y$ run over the elements of an admissible
affinoid cover of~$\W$, the corresponding $Z_{\phi,Y}$ glue together
to give the \emph{spectral curve}
$Z_\phi\subset\W\times\A^1$ associated to $\phi$.

Going back to $Y\to\W$ arbitrary, for every factorization $P_{\phi,Y}=Q(T)S(T)$
where $Q=1+\ldots$ is a polynomial of degree~$n$ with leading term
a unit, and such that $Q$ and $S$ are relatively prime, there is a
decomposition of $Z_{\phi,Y}$ into two disjoint subspaces, the one
corresponding to~$Q$ being finite and flat of degree~$n$. The
submodule $N$ of $M_Y[r]$ corresponding to $Q(T)$ via Coleman's Riesz
theory (Theorem~A4.3 of~\cite{Coleman})
is a free $R$-module of rank~$n$, and is stable under all the Hecke operators. 
The Hecke operators acting on this subspace generate a finite
free $R$-algebra which is hence an affinoid algebra,
and the associated rigid space naturally lives
over $Z_{\phi,Y}$. The eigencurve $\E$ is built by
glueing all such spaces together, as $Y$ ranges over admissible
affinoid opens in $\W$. The key technical difficulty in this construction
is verifying that the induced cover of $Z_\phi$ is admissible,
and this problem was solved in Proposition~A5.8 of~\cite{Coleman}.

\begin{lemma}
\label{boundedhecke}
Let $Y$ be an affinoid subdomain of $\W$ and choose $r>0$
sufficiently small so that the Hecke operator $T_\ell$ is a well-defined
endomorphism of $M_Y[r]$. Then, possibly after shrinking~$r$ again,
we have $||T_\ell||\leq 1$.
\end{lemma}
\begin{proof} This comes from an explicit analysis of the formula
used to describe $T_\ell$. Note that, unlike the classical case,
$T_\ell$ is not defined as a correspondence in general weight~$\kappa$,
because of the slightly unnatural definition of an overconvergent
eigenform of weight $\kappa$. On the other hand, the definition
as a correspondence does work well in weight~0, and one can deduce
from this that $\sum a_nq^n\mapsto\sum a_{n\ell}q^n$
is a well-defined map from $r$-overconvergent functions of weight~0
and level~1 to level $r$-overconvergent functions of weight~0 and level~$\ell$,
and furthermore that this map has norm at most~1.
Now the lemma follows from the explicit definition of $T_\ell$
at weight $\kappa$ given on p463 of~\cite{Coleman},
with the proviso that this definition only works near the centre
of weight space, so every occurrence of $E^s$ should be replaced by
the Eisenstein family $E_\kappa$, and every occurrence of
$e_\ell^s$ should replaced by $E_\kappa(q)/E_\kappa(q^\ell)$,
a function which is proved to be overconvergent in Proposition~2.2.7
of~\cite{Eigencurve}, and which has $q$-expansion congruent to~1
modulo the maximal ideal of $\Ot_2$. The reason one might have to
shrink $r$ again is that we need to guarantee that the supremum
norm of $E_\kappa(q)/E_\kappa(q^\ell)$ is at most~1 on
$X_0(\ell)_{\geq 2^{-r}}$.
\end{proof}
\begin{lemma}\label{1wt0}
If $\kappa\not=0$ then $E_\kappa$ is not the $q$-expansion
of a function on $X[0]$. Equivalently, $1$ is not the $q$-expansion
of a 2-adic modular form of weight $\kappa$ for any non-zero $\kappa$.
\end{lemma}
\begin{proof} For $\kappa$ a positive even integer this follows
from Corollary~4.5.2 of~\cite{Katz}. We reduce to this case. Assume
$E_\kappa$ is a function on $X[0]$. Then theorem~2.2.2 of~\cite{Eigencurve}
implies that $E_{\kappa^2}$ is as well, and hence we may assume
that $\kappa=\langle\cdot\rangle^{k'}$ with $k'\in\C_2$ and furthermore
we may assume that $k'$ is sufficiently close to zero to ensure
that $E_{k'}-1$ has $q$-expansion divisible by, say, 16 in $\Ot_2[[q]]$.
Now choose $k\in\Z_{>0}$ with $k=k'\alpha$ and $|\alpha|<1$, and
consider the function $(E_{k'})^{\alpha}$ on $X[0]$. By Corollary~B4.5.2
of~\cite{Coleman} we see
that $E_k/(E_{k'})^{\alpha}$ is an overconvergent function on $X[0]$,
and this reduces us to the case we have dealt with already.
\end{proof}
Fix a weight $\kappa$, and for a Hecke operator
$T\in\{U,T_3,T_5,T_7,T_{11},\ldots,T_\ell,\ldots\}$ define $\lambda(T)$ to be
the eigenvalue of $T$ acting on $E_\kappa$ (so $\lambda(U)=1$
and $\lambda(T_\ell)=1+\kappa(\ell)/\ell$.)
\begin{cor}\label{noneis}
If $f$ is a non-zero
cuspidal eigenform of weight $\kappa$, then there is a Hecke
operator $T\in\{U,T_3,T_5,T_7,T_{11},\ldots\}$ such that
$Tf\not=\lambda(T)f$.
\end{cor}
\begin{proof}
By standard results on how Hecke operators act on
$q$-expansions, we see that any counterexample to the lemma
must be of the form $f=c(q+q^2+\lambda(T_3)q^3+\cdots)$ with $c\not=0$.
If $\kappa\not=0$
then this eigenform is of the form $c'(E_\kappa-1)$,
and hence the $q$-expansion 1
is in the linear span of $E_\kappa$ and $f$, and we deduce that~1
is an overconvergent modular form of weight~$\kappa$, contradicting
Lemma~\ref{1wt0}. It remains to deal with the case $\kappa=0$.
Yet, as noted in Lemma~4 of~\cite{CGJ}, a result of Serre implies that
the form $f$ is not even a $p$-adic modular form, and thus certainly
not an overconvergent eigenform.
\end{proof}

\begin{lemma}\label{ce} The 2-adic level~1 eigencurve
$\E$ can be written as a disjoint union $\Ee\coprod\Ec$,
with $\Ee$, the Eisenstein component, mapping isomorphically
down to $\W$ via the projection, and $\Ec$ being the eigencurve
constructed from spaces of cuspidal overconvergent modular forms
via the argument above.
\end{lemma}
\begin{proof} This is no doubt well-known but we write down a proof
for lack of a reference. Let $P(T)$ denote the characteristic polynomial
of $U$ in $\Ot(\W)[[T]]$. Now $P(1)=0$ because it is a function
on weight space that vanishes at all classical even weights $k\geq2$,
which are Zariski-dense in $\W$ (it vanishes because the Eisenstein
series is an eigenform with eigenvalue~1). Write $P(T)=(1-T)P^0(T)$.
Set $z:=P^0(1)$. Now $z$ is not identically zero, because if it were
then there would be a cuspidal overconvergent eigenform of
weight~4 with $U$-eigenvalue~1 and such a thing would be
classical. However the level~2 weight~4 Eisenstein series which vanishes at
infinity does not have the right eigenvalue, and neither do any cusp
forms because this would contradict the Weil bounds.
So the zeroes of $z$ form a Zariski-closed
subset of weight space which is not all of $\W$.
Let $\W^\times$ denote the complement of this
set, so $\W^\times$ is open and dense in $\W$.
Over $\W^\times$ we know that $(1-T)$ is coprime to $P^0(T)$,
and one deduces that the spectral curve $Z_U$
is the disjoint union of the component $Z_U^{\mathrm{Eis}}$
corresponding to the $U$-eigenvalue~1,
and its complement, corresponding to cusp forms.
Moreover the construction of the eigencurve over the spectral
curve gives, over $Z_U^{\mathrm{Eis}}$, a component of the eigencurve
isomorphic to $\W^\times$, since the associated Hecke algebra
is of rank~1. Hence over $\W^\times$ the eigencurve is a disjoint
union of a component $\Ee$ isomorphic to $\W^\times$ and its
complement, $\Ec$.

We must extend this construction now to $\W$. We remark
that in the case we are interested in it is almost certainly
the case that $\W^\times=\W$, and this would follow from the well-known
fact that Hida theory and Coleman theory are compatible;
unfortunately we have been unable to find an explicit reference for this
that applies for small primes or for
weights that are not in $\Hom(\Z_2^\times,\Z_2^\times)$, so we give
a self-contained proof. The trick is to change our choice
of compact
operator. If $\kappa\in\W\backslash\W^\times$ then there is a cuspidal
eigenform of weight $\kappa$ with $U$-eigenvalue~1, but we shall
construct another compact operator $U'$ such that the eigencurves
constructed via $U$ and $U'$ are isomorphic, $E_\kappa$ is an
eigenvector for $U'$ with eigenvalue $\alpha$,
and furthermore $(1-\alpha T)$ divides the characteristic power series 
of $U'$ on $M_\kappa[r]$ precisely once. The existence of such
a $U'$ implies that $Z_{U'}$ splits up as the disjoint union of an
Eisenstein component and a cuspidal component over a neighbourood of $\kappa$,
and hence $\E$ also splits up as a disjoint union of $\Ec$
and $\Ee$ over this neighbourhood, which is what we need
to finish the proof.

It remains to construct such a $U'$.
Chose $\kappa$ in $\W\backslash\W^\times$
and consider the space $V$ of overconvergent cusp forms of weight $\kappa$
annihiliated by $(U-1)$. This space is finite-dimensional
and non-zero. Furthermore, by Corollary~\ref{noneis},
for any $0\not=v\in V$ there exists
an odd prime $\ell$ such that $T_\ell v\not=(1+\kappa(\ell)/\ell)v$.
Choose a basis $\{e_1,e_2,\ldots,e_n\}$
of $V$ such that all the $T_\ell$ are in upper
triangular form, and for each $e_i$ choose a Hecke operator $T_{\ell_i}$
such that $T_{\ell_i}e_i\not=(1+\kappa(\ell_i)/\ell_i)e_i$. It is easy
now to find a linear combination $T_0:=\sum_ic_iT_{\ell_i}$
of these Hecke operators
such that if $T_0E_\kappa=\lambda E_\kappa$ then $\lambda$ is
not an eigenvalue of $T_0$ on $V$. For $N$ sufficiently large
we have $||p^NT_0||<1$ and hence $1+p^NT_0$ is invertible on
$M_\kappa[r]$. We claim that for some such $N$ the Hecke operator
$U':=U(1+p^NT_0)$ suffices. The eigencurves
constructed using $U$ and $U'$ are isomorphic above a small neighbourhood
of $\kappa$ in $\W$, by the arguments of Corollary~7.3.7 of~\cite{Eigencurve}
(applied to the neighbourhood of $\kappa$ rather than all of weight
space, and noting that the argument
does not rely on any of the deformation theory of Galois
representations presented earlier in~\cite{Eigencurve} and hence does not need
the assumptions $N=1$ and $p>2$.) It remains to check that we can choose
$N$ such that if $U'E_\kappa=\alpha E_\kappa$ then the generalised
$\alpha$-eigenspace for $U'$ on $M_\kappa[r]$
is precisely 1-dimensional (and hence spanned by $E_\kappa$). It suffices
to verify this on the $U$-ordinary subspace of $M_\kappa[r]$ (which is
equal to the $U'$-ordinary subspace of $M_\kappa[r]$), as
$\alpha$ is a unit.
The ordinary subspace splits as a direct
sum of the Eisenstein subspace and the cuspidal part,
which in turns splits into the sum of
the generalised $U$-eigenspace $V_1$ where the $U$-eigenvalue is~1,
and the direct sum $V_2$ of the other generalised $U$-eigenspaces.
On $V_1$ we have to verify that no eigenvalue of $UT$ is $\alpha$,
which follows without too much trouble from our construction of $T$,
whatever the value of~$N$. Finally the space $V_2$ is finite-dimensional
and 1 is not an eigenvalue of $U$ on this space. On the other hand,
as $N$ tends to infinity we see that $\alpha$ tends to~1 and $U'$
tends to~$U$, so for $N$ large enough there will also be no
eigenvalues equal to $\alpha$ on $V_2$. This completes the proof.
\end{proof}
Finally we need a result that says that $\Ec$ represents a functor
on rigid spaces over weight space. Again this result seems to be
known but we know of no reference. If $\W_i$ is an affinoid subdomain
of $\W$ then we let $S_{\W_i}[r]$ denote the $r$-overconvergent cusp
forms of weight $\W_i$. Now let $Y$ denote any rigid space
over $\W$. We say that $F=\sum_{n\geq1}a_nq^n\in\Ot(Y)[[q]]$ is a
\emph{normalised overconvergent finite slope cuspidal eigenform of weight~$Y$}
if $F=q+O(q^2)$, if $a_2\in\Ot(Y)^\times$,
and furthermore if we can write $Y$ as an admissible
union of affinoids $Y_i$ such that for each $i$ there exists
$r_i>0$ and an affinoid subdomain $\W_i$ of $\W$ with $Y_i\to\W_i$,
such that $F$ is the $q$-expansion of an element
in $S_{Y_i}[r]:=S_{\W_i}[r]\widehat\otimes_{\Ot(\W_i)}\Ot(Y_i)$.
Let $\F(Y)$ denote the functor on rigid spaces over $\W$, sending $Y/\W$
to the set of normalised overconvergent finite slope cuspidal eigenforms of
weight~$Y$.
\begin{lemma}
\label{rep}
This functor is represented by $\Ec$.
\end{lemma}
\begin{proof} We need to exhibit functorial bijections
$\Ec(Y)=\F(Y)$ for all $Y$, which we do by writing down
canonical maps in both directions. Let us first start with
a map $\beta:Y\to\Ec$ and concoct a finite slope cuspidal
overconvergent eigenform. Recall that $\Ec$
is equipped with functions $T_1$, $T_2$, $T_3,\ldots$
and given $\beta:Y\to\Ec$ we define $a_n=\beta^*(T_n)$
and set $F=\sum_{n\geq1}a_nq^n$. We claim that this
is indeed a normalised
overconvergent finite slope cuspidal eigenform. It suffices to
verify this on an admissible affinoid cover of~$Y$, and hence we
may assume that $Y=\Sp(A)$ is affinoid and the map $Y\to\Ec$ has
image in $\Sp(\T)$, where $\T$ is one of the Hecke algebras
used to define the eigencurve via the spectral curve $Z_U$.
Now $\T$ is the Hecke algebra
corresponding to a finite rank space of overconvergent
cuspidal modular forms $M$ over an affinoid~$R$, and Coleman
proves on p465 of~\cite{Eigencurve} 
that the usual $R$-linear pairing $\T\times M\to R$
defined by $(t,m)\mapsto a_1(tm)$
is perfect. Because $\T$ and $M$ are free $R$-modules of finite rank,
this pairing remains perfect when one tensors up to $A$,
and we deduce that the map $\T\to A$ of $R$-modules corresponds
canonically to an element of $M\otimes_RA=M\widehat\otimes_R A$ with
$q$-expansion $\sum_{n\geq 1} a_n q^n$, and in particular to
a normalised cuspidal overconvergent eigenform. It suffices to
prove that this eigenform has finite slope, but this is clear
because, by definition, the characteristic polynomial of $U$
on~$M$ has constant coefficient equal to a unit in $R$, and hence $U=T_2$
is invertible, thus $a_2$ is also invertible.

The construction the other way is just a case of ensuring that
the argument above can be reversed. If $Y/\W$ is a rigid space
and $F=\sum a_nq^n$ is a normalised cuspidal finite slope overconvergent
eigenform over $Y$ then we must show that there is a unique
map $Y\to\Ec$ such that $a_n$ is the pullback of $T_n$ for all $n$.
Again it suffices to do
this on an admissible affinoid cover of $Y$ so again we may assume $Y=\Sp(A)$
is affinoid, that the map $Y\to\W$ has image in an affinoid $\W_i$,
and that $\sum a_nq^n$ is an element of a space
$S_Y[r]$ of $r$-overconvergent cusp forms of weight~$Y$. 
Let $P_{U,Y}(T)$ denote the characteristic power series of $U$
on $S_Y[r]$. The factor $(1-a_2T)$ of $P_{U,Y}(T)$, corresponding
to our finite slope eigenform, cuts out a closed subspace $Z_F$ of $Z_{U,Y}$,
which maps down isomorphically onto $Y$ under the canonical projection
$Z_{U,Y}\to Y$, as $a_2\in\Ot(Y)^\times$. Note that $Z_F$ may not
be disconnected from the closure of its complement in $Z_{U,Y}$.
The admissible cover of $Z_U$ in Proposition~A5.8 of~\cite{Coleman}
pulls back to an admissible
cover of $Z_{U,Y}$ and hence to an admissible cover of $Z_\F$ and thus
of $Y$. Replacing $Y$ by an element of this admissible cover, 
we may assume that there exists a factorization $P_{U,Y}(T)=Q(T)S(T)$
with $(Q(T),S(T))=1$, $Q(T)=1+O(T)$ a polynomial with leading term a unit, and
$(1-a_2T)|Q(T)$. This factorization induces a $U$-invariant
decomposition $S_Y[r]=N\oplus E$ with $N$ free of finite rank
over $\Ot(Y)$.
We may write $F=F_N+F_E$ via this decomposition, and both $F_N$
and $F_E$ will be eigenvectors for $U$ with eigenvalue $a_2$. But
$S(T)$, the characteristic power series of $U$ on $N$, is coprime
to $Q(T)$ and hence to $1-a_2T$, so $F_E=0$. We deduce that $F\in N$.
Now if $\T$ denotes the Hecke algebra associated to $N$ then $F\in N$
induces an $\Ot(Y)$-algebra homomorphism $\T\to\Ot(Y)$, and it is a standard
calculation, using the fact that $F=q+\cdots$ is an eigenform, that
this map is in fact a ring homomorphism. This ring
homomorphism induces a map $Y\to\Sp(\T)$ and hence $Y\to\Ec$.

Finally, it is elementary to verify that both constructions are inverse
to one another.
\end{proof}
\begin{remark}\label{76}
 Lemmas~\ref{boundedhecke} and~\ref{1wt0} are true for general Coleman-Mazur
eigencurves (with the same proofs!).
Corollary~\ref{noneis} is true for regular primes
but will not be true at weights corresponding to zeros of the $p$-adic
$L$-function, because the corresponding Eisenstein series is cuspidal.
Similarly for Lemma~\ref{ce} --- the proof works for
regular primes but the cuspidal and Eisenstein components of the
eigencurve will meet for irregular primes, as can be seen from
the main theorem of~\cite{Emerton} and the well-known compatibility
of Hida theory and Coleman theory. On the other hand Lemma~\ref{rep} is
true for
general Coleman-Mazur eigencurves --- one can define $\Ec$ using families of
cuspidal overconvergent modular forms, rather than as a component of $\E$.
\end{remark}

\section{There are not too many holes in the eigencurve.}

We begin with a simple rigid-analytic lemma that forms the basis
to our approach.
Let~$X$ be a connected affinoid variety, and let~$V$
be a non-empty admissible open affinoid subdomain of~$X$.
Let~$B=\Sp(\C_2\langle T\rangle)$ denote the closed unit disc,
and let $A=\Sp(\C_2\langle T,T^{-1}\rangle)$ denote its ``boundary'',
the closed annulus with inner and outer radii both~1.

\begin{lemma} If $f$ is a function on $V\times B$ and the restriction
of $f$ to $V\times A$ extends to a function on $X\times A$,
then $f$ extends to a function on $X\times B$.
\label{lemma:rigid}
\end{lemma}
\begin{proof} We have an inclusion $\Ot(X)\subseteq\Ot(V)$, as~$X$
is connected,
and we know $f\in\Ot(V)\langle T\rangle$ and $f\in\Ot(X)\langle T,T^{-1}\rangle$.
But the intersection of these two rings is $\Ot(V)\langle T\rangle$.
\end{proof}
Let $\E$ denote the 2-adic eigencurve of tame level~1, and let $\W$
denote 2-adic weight space. Let $B$ denote the closed unit disc
and let $B^\times$ denote $B$ with the origin removed.
Suppose we have a map $\phi:B^\times\to\E$ such that the induced
map $B^\times\to\W$ extends (necessarily uniquely) to a map $B\to\W$.
Let $\kappa_0\in\W(\C_2)$ denote the image $0\in B(\C_2)$ under this map.
The theorem we prove in this section is
\begin{theorem}\label{main:one} If
$\kappa_0\notin\{\langle\cdot\rangle^{-12s}:2s\in\Z_2^\times\}$
then the map $\phi:B^\times\to\E$ extends to a map $B\to\E$.
\end{theorem}
\begin{proof}
Recall from Lemma~\ref{ce}
that $\E=\Ec\coprod\Ee$. If the image of $\phi$ is contained
in $\Ee$ then the theorem is automatic, since the projection
$\Ee\to\W$ is an isomorphism. Hence we may assume that $\phi:B^\times\to\Ec$.
If $|\kappa_0(5)-1|>1/8$ then we are finished by the main theorem
of~\cite{buzzkill}. Assume from now on that $|\kappa_0(5)-1|\leq 1/8$.
By Lemma~\ref{rep} the map $\phi$ gives rise to 
a family $\sum a_n q^n$ of overconvergent
eigenforms over $B^\times$. By Lemma~\ref{boundedhecke} the supremum
norm of each $a_n$ is at most~1 and, analogous to the analysis of
isolated singularities of holomorphic functions, one checks easily
that this is enough to ensure that each $a_n$ extends to a function on
$B$.
We now have a formal power
series $\sum_{n\geq1}a_nq^n$ in $\Ot(B)[[q]]$. We
next claim that this formal power series is an overconvergent form
of weight $B$ --- indeed, it is not too difficult to establish how
overconvergent it is. We are assuming $|\kappa_0(5)-1|\leq1/8$
and hence $\kappa_0=\langle\cdot\rangle^{-12s}$ with $|s|<4$.
Now assume also that $2s\not\in\Z_2^\times$. Set $r=\frac{3+\nu(2s)}{12}$.
After shrinking $B$ if necessary, we may assume that 
for all $b\in B$
we have $\kappa_b=\langle\cdot\rangle^{-12s'}$ with $|s-s'|\leq 1$.
By Lemma~\ref{lemma:smalldisc} we have $(\kappa_b,r)\in\X$ for all $b\in B$,
and by Corollary~\ref{cor:convergesfam} we see that on the boundary
of $B$ our function $\sum a_nq^n$ is $r$-overconvergent, it being
a finite slope eigenform for $U$ here. By Lemma~\ref{boundedhecke}
the coefficients $a_n$ are all bounded by~1 on all of $B$.
Now applying Lemma~\ref{lemma:rigid}
with $X=X[r]$ and $V$ a small disc near infinity such such that $q$
(the $q$-expansion parameter) is a well-defined function on $V$,
we deduce that $\sum a_nq^n$ is $r$-overconvergent on all of $B$.

Next we show that $a_2\in\Ot(B)^\times$. It suffices to prove
that $a_2(0)\not=0$, as we know that $a_2(b)\not=0$ for all $0\not=b\in B$.
But $\sum a_n(0)q^n=q+\ldots$ is an $r$-overconvergent form of
weight $\kappa_0$, so by Theorem~\ref{52} we deduce $a_2(0)\not=0$.
Hence $a_2\in\Ot(B)^\times$ and $\sum a_nq^n$ is an overconvergent
cuspidal finite slope eigenform of weight~$B$. We finish the proof
by applying Lemma~\ref{rep} once more, giving us a map $B\to\Ec$.
\end{proof}

% Note to Frank: I removed the table. The reason is that ``I found a formula
% for the entries'': the finite slopes ones go up to 1/2+mu(s)/6-epsilon,
% and the infinite slopes ones don't go up to 1/4+nu(2s)/12. Lemma 5.2
% then tells us what we need.
%
%\begin{table}
%\begin{center}
%\begin{tabular}{|c|c|c|}
%\hline
%$\mu$ & Infinite slope Eigenforms & 
%Finite slope eigenforms  \\
%\hline
%$-2 < \mu < -1$ & $1/6 + (\mu+2)/12$ & $1/6 + (\mu+2)/6$ \\
%\hline
%$-1 \le \mu$, $2s \notin \Z_2$ & $1/3 - \epsilon$ &  $1/3 + (\mu + 1)/6$\\
%\hline
%$2s \in \Z_2$, $\mu = -1$, $6s \notin -\N$ & $1/3$ & $1/3$\\
%\hline
%$2s \in \Z_2$, $\mu > -1$ & $1/3$ & $1/3 + (\mu + 1)/6$\\
%\hline
%%$6s \in 2 \Z$, $\mu = 0$ & $1/3$ & $1/2$\\
%%\hline
%\end{tabular}
%\end{center}
%\caption{Radius of convergence for eigenforms of finite and 
%infinite slope.}
%\label{tab:mytable}
%\end{table}

\section{There are no holes in the eigencurve}

In the previous section we showed that if there are any holes in the
eigencurve, then they lie above weights of the form
$\{\langle\cdot\rangle^{-12s}:2s\in\Z_2^\times\}$.
To show that in fact there are no holes in the eigencurve,
we redo our entire argument with a second, even more non-standard,
twist and show that using this twist we may deduce that
the only holes in the eigencurve lie above the
set $\{\langle\cdot\rangle^{2-12s}:2s\in\Z_2^\times\}$.
Because there is no $s\in\frac{1}{2}\Z_2^\times$ such that
$\frac{12s-2}{12}\in\frac{1}{2}\Z_2^\times$ this finishes the
argument. We sketch the details.

Let $E_2=1 + 24q + 24q^2 + 96q^3 +\ldots$ denote the holomorphic
Eisenstein series of weight~2 and level $\Gamma_0(2)$.
We define $\X'=\{(\kappa\langle\cdot\rangle^2,r):(\kappa,r)\in\X\}$.
If $|s|<8$ then set $\kappa'=\langle\cdot\rangle^{2-12s}$. 
If $r$ is such that $(\kappa',r)\in\X'$, we define $M'_{\kappa'}[r]$ to be
the vector space of formal $q$-expansions $F\in\C_2[[q]]$ such
that $Fh^s/E_2$ is the $q$-expansion of an element of $M_0[r]$.
For $r>0$ sufficiently small this definition is easily checked
to coincide with the usual definition. We shall be using this definition
with $r$ quite large and again we neglect to verify whether the
two definitions coincide in the generality in which we use them.
We give $M'_{\kappa'}[r]$ the Banach space structure such that
multiplication by $h^s/E_2$ is an isometric isomorphism
$M'_{\kappa'}[r]\to M_0[r]$,
and endow $M'_{\kappa'}[r]$ once and for all with the orthonormal basis
$\{E_2h^{-s},E_2h^{-s}(2^{12r}f),E_2h^{-s}(2^{12r}f)^2,\ldots\}$.
Note that the reason that this definition gives us more than
our original definition of $M_\kappa[r]$
is that if $k$ is an even integer with $2||k$ then $(k,1/3)\not\in\X$
but $(k,1/2-\epsilon)\in\X'$, so we can ``overconverge further''
for such weights.

If $\theta=q(d/dq)$ is the operator on formal $q$-expansions,
then one checks that $U \theta = 2 \theta U$. Moreover,
it is well-known that $\theta f=fE_2$ and hence
$\theta f^j = j f^j E_2$ for any $j\geq0$.
Hence our formulae for the coefficients of $U$ acting on $M_0[r]$
will give rise to formulae for the coefficients
of $U$ acting on $M'_2[r]$, which was the starting point
for the arguments in section~3. We give some of the details
of how the arguments should be modified.
If $m\in\Z_{\geq0}$ and $k'=2-12m$
then we define a continuous operator $U'_{k'}$ on $M_0[r]$ by
$U'_{k'}(\phi)=E_2^{-1}h^mU(E_2h^{-m}\phi)$. One checks
that this is indeed a continuous operator by verifying that it
has a basis $(u'_{ij}(m))_{i,j\geq0}$ defined by $u'_{ij}(m)=0$
for $2i<j$ or $2j-i+3m<0$, $u'_{00}(0)=1$, and
$$u'_{ij}(m) = \frac{3(i+j + 3m - 1)! (i + m) 2^{8i - 4j + 12r(j-i)}}
{(2i-j)! (2j-i + 3m)!}$$
otherwise.
One checks that for $i,j$ fixed there is a polynomial
$u'_{ij}(S)$ interpolating $u'_{ij}(m)$ and that for $|s|<8$
with $\mu=\min\{v(s),0\}$ we have
$v(u'_{ij}(s))\geq(\mu+3-6r)(2i-j)+6rj$ as before.
Hence for $|s|<8$, $\kappa'=\langle\cdot\rangle^{2-12s}$
and $r\in\Q$ such that $(\kappa',r)\in\X'$, the matrix $(u'_{ij}(s))_{i,j\geq0}$
defines a compact operator $U'(s)$ on $M_0[r]$. Furthermore we have
$U'(s)(\phi)=E_2^{-1}h^sU(E_2h^{-s}\phi)$, and in particular
$U:M'_{\kappa'}[r]\to M'_{\kappa'}[r]$ is well-defined and compact. Moreover
$U'(s)$ increases overconvergence and any eigenvector for $U'(s)$
on $M_0[r]$ with non-zero eigenvalue extends to $M_0[r']$ for
any $r'$ such that $0<r'<1/2+\mu(s)/6$. Finally, these arguments
also work for families of modular forms and the analogue of
Corollary~\ref{cor:convergesfam} remains true in this setting.

Similar arguments work in section 4. One checks
that $2V\theta=\theta V$ and hence
$V U \theta = 2 V  \theta U = \theta V U$. Hence $\theta$
commutes with $W$ and one now deduces from our explicit formulae
for $W$ in weight $-12m$ that in weight $2-12m$ the matrix
for $W$ is given by $W_k = [\eta'_{ij}]$, where:
$$\eta'_{ij} = \frac{(2i + j-1 + 6m)! 3(i+2m) \cdot 2^{(4 -12 r)(i-j)} (-1)^i}
{(i-j)! (i+2j+6m)!}.$$
We remark that the only difference in this formula is that $(j + 2m)$
has been replaced by $(i + 2m)$. One finds that the arguments
at the end of this section apply
\emph{mutatis mutandis} in this case.

The analogue of Theorem~\ref{52} is that if $|s|<4$ and $2s\not\in\Z_2^\times$
and $\kappa'=\langle\cdot\rangle^{2-12s}$ then an overconvergent infinite
slope form of weight $\kappa'$ is not $r$-overconvergent, for
$r=\frac{3+\nu(2s)}{12}$. The proof follows the same strategy, although
some of the lemmas in section~6 need minor modifications; for example
in Lemma~\ref{combinatorics} we set $x_1=\frac{i+2j+6s}{3(i+2s)}\eta'_{ij}$
and $x_2:=\frac{2i+j+6s}{3(i+2s)}\eta'_{ij}$, and the result
follows as $\eta'_{ij}=2x_2-x_1$. Note that $E_2$ can be regarded
as an element of $H^0(Y[1/3],\omega^{\otimes 2})$ so that 
Lemma~\ref{integercase} does not need modification.

We deduce our main theorem:
\begin{theorem} If $\phi:B^\times\to\E$ and the induced map $B^\times\to\W$
extends to a map $\psi:B\to\W$, then $\phi$ extends to a map $B\to\E$.
\end{theorem}
\begin{proof} If
$\psi(0)\not\in\{\langle\cdot\rangle^{-12s}:2s\in\Z_2^\times\}$ then
we use Theorem~\ref{main:one}, and if it is then we use the
modification explained above.
\end{proof}

\noindent \it Email addresses\rm:\tt \  buzzard@imperial.ac.uk
\hskip 22mm \tt \ fcale@math.harvard.edu

\end{document}